\documentstyle[11pt]{amsart}
\catcode`\@=11

\long\def\@savemarbox#1#2{\global\setbox#1\vtop{\hsize\marginparwidth 
  \@parboxrestore\tiny\raggedright #2}}
\newcommand\lref[1]{\ref{#1}%
\@ifundefined{r@DisplaY #1}{}{(#1)}}% Prints label as well as
\newcommand\fakelabel[2]{\@bsphack\if@filesw {\let\thepage\relax
   \newcommand\protect{\noexpand\noexpand\noexpand}%
\xdef\@gtempa{\write\@auxout{\string
      \newlabel{#1}{{#2}{\thepage}}}}}\@gtempa
   \if@nobreak \ifvmode\nobreak\fi\fi\fi\@esphack}

\catcode`\@=12

\def\Empty{}
\newcommand\oplabel[1]{
  \def\OpArg{#1} \ifx \OpArg\Empty {} \else
  	\label{#1}
  \fi}
\newtheorem{theoremSt}{Theorem}[section]
\newtheorem{corollarySt}[theoremSt]{Corollary}

\newtheorem{exampleSt}[theoremSt]{Example}
\newtheorem{exerciseSt}[theoremSt]{Exercise}

\newcommand\MakeStEnv[1]{
  \newenvironment{#1}[1]{%    environment without explicit label
  \begin{#1St} \oplabel{##1}%
  \global\def\CrntSt{\thetheoremSt}%
}{ 
  \end{#1St} }
  \newenvironment{#1+}[1]{%   environment with explicit label
  \begin{#1St} \label{##1}%
  \label{DisplaY ##1}%
  \global\def\CrntSt{\thetheoremSt}%
  \def\Labl{##1}\ifx\Labl\Empty{} \else {\em (\Labl)\,}\fi%
}{ 
  \end{#1St} }
}
\MakeStEnv{theorem}
\MakeStEnv{corollary}
\MakeStEnv{proposition}
\MakeStEnv{lemma}
\MakeStEnv{define}
\MakeStEnv{conjecture}

\newlength{\saveu}

\newcommand{\finishproof}[1]{ 
  \def\FPArg{#1}
  \ifx\FPArg\Empty
  	\newcommand\FPArg{\CrntSt}  \fi
  \smallbreak\noindent\makebox[\textwidth]{\hfill\fbox{\FPArg}}
  \medbreak\noindent
}

\newcommand\DD{{\cal D}}

\newcommand\FF{{\cal F}}

\newcommand\HH{{\cal H}}

\newcommand\LL{{\cal L}}

\newcommand\bbR{{\mathord{\text{I\kern-2pt R}}}}        % Fake blackboard
\newcommand\bbH{{\mathord{\text{I\kern-2pt H}}}}        % Fake blackboard

\newcommand\C{{\bold C}}

\newcommand\R{{\bold R}}

\newcommand\bigrightarrow[1]{\hbox to #1{\rightarrowfill}}
\newcommand\bigleftarrow[1]{\hbox to #1{\leftarrowfill}}

\newcommand\semidir{\mathrel{\hbox{\vrule depth-.03ex 
	height1.1ex\kern-0.15em$\times$}}}

\newcommand{\pair}[1]{\langle #1\rangle}

\numberwithin{equation}{section}

\renewcommand{\Re}{\mathop{\hbox{\rm Re}}}
\begin{document}

\title[]{Quasi-rigidity of Hyperbolic $3$-Manifolds and Scattering Theory} 
\author{David Borthwick}
\address{Department of Mathematics, University of Michigan,
Ann Arbor, MI 48109}
\author{Alan McRae}
\address{Department of Mathematics, University of Indiana, Bloomington, 
	IN  47405}
\author{Edward C. Taylor}
\address{Department of Mathematics, University of Michigan,
Ann Arbor, MI 48109}

\date{\today}

\thanks{First author supported in part by NSF grant DMS-9401807.
Third author partially supported by a University of Michigan Rackham
Fellowship.}

\begin{abstract} 

Suppose $\Gamma_{1}$ and $\Gamma_{2}$ are two convex co-compact
co-infinite volume discrete subgroups of $PSL(2,{\bf C})$, so that there
exists a $C^{\infty}$ diffeomorphism 
$\psi: \Omega(\Gamma_{1})/\Gamma_{1} \rightarrow
\Omega(\Gamma_{2})/\Gamma_{2}$ that induces an isomorphism $\phi:
\Gamma_{1} \rightarrow \Gamma_{2}.$  For fixed  $s \in {\bf C}$, let
$S_{i}(s)$ be the scattering operator on $\Omega(\Gamma_{i})/\Gamma_{i}\;
(i=1,2)$.   Define $S_{rel}(s) = S_{1}(s) - \psi^{*}S_{2}(s)$, where
$\psi^{*}S_{2}(s)$ is the pull-back of $S_{2}(s)$ to an operator acting on
the appropriate complex line bundle over $\Omega(\Gamma_{1})/\Gamma_{1}.$  
Our main result is: if the operator norm of $S_{rel}(s)$ is
$\epsilon$-small, then the $\Gamma_{i}$ are
$K(\epsilon)$-quasi-conformally conjugate and the dilatation  $K(\epsilon)$
decreases to $1$ as $\epsilon$ decreases to $0$. 
  
\end{abstract}

\maketitle

\section{Statement of Results}

Geometrically finite Kleinian groups uniformizing infinite volume
hyperbolic $3$-manifolds exhibit a 
rich deformation theory  due to work of Ahlfors, Bers, Kra, Marden, Maskit,
Thurston and others.  The purpose of this note is to introduce 
scattering theory as an analytic tool in the study of the deformations of
complete geometrically finite hyperbolic structures. 
The results in this paper are restricted to a sub-class of geometrically
finite Kleinian groups called {\em convex co-compact groups}, i.e. those
containing no parabolic subgroups. 

Assume that $\Gamma$ is a convex co-compact, torsion-free Kleinian
group with non-empty regular set $\Omega(\Gamma)$ (see Section 2 for
definitions).  The compact  (possibly disconnected) 
quotient surface $\Omega(\Gamma)/\Gamma$ is the conformal boundary at
infinity of the hyperbolic 3-manifold $M(\Gamma) = {\bf H}^3/\Gamma$.  To
the Laplacian 
$\Delta$ on $M(\Gamma)$ we can associate a scattering operator acting
on sections of certain complex line bundles over the conformal boundary 
$\Omega(\Gamma)/\Gamma$.  These sections are most conveniently described as
automorphic forms on $\Omega(\Gamma)$.  For a complex parameter $s$,   
let $\FF_{s}(\Gamma)$ be the space of automorphic forms of weight $s$ on
$\Omega(\Gamma)$ (see Section 3 for the definition).
The scattering operator $S(s)$ is a pseudodifferential operator
\cite{Man89}  with known singularity  
mapping $\FF_{2-s}(\Gamma) \to \FF_{s}(\Gamma)$.  
For $\Re s = 1$ we have a natural $L^2$ inner product on $\FF_s(\Gamma)$,
so we can complete these spaces to form Hilbert spaces.

Now take two convex co-compact, torsion-free Kleinian groups
$\Gamma_{i=1,2}$, with $\Omega(\Gamma_{i=1,2}) \neq \emptyset$.  Assume
there exists an orientation-preserving diffeomorphism $\psi :
\Omega(\Gamma_{1}) \rightarrow \Omega(\Gamma_{2})$ that induces an
isomorphism $\phi: \Gamma_{1} \rightarrow \Gamma_{2}$. 
Denote by $S_{2}(s)$  the scattering operator
acting on sections over $\Omega(\Gamma_{2})/\Gamma_{2}$.  Recall the
diffeomorphism $\psi$ descends to a diffeomorphism $\psi:
\Omega(\Gamma_{1})/\Gamma_{1} \rightarrow \Omega(\Gamma_{2})/\Gamma_{2}$.
Thus we can pull  $S_{2}(s)$ back, via the diffeomorphism $\psi$, to an
operator  taking acting on sections over $\Omega(\Gamma_{1})/\Gamma_{1}$. 
Denote this pull-back of $S_{2}$ by $\psi^{*}S_{2}(s)$.  
Perry \cite{Per95} shows that if for some $\Re s = 1, s\ne 1$, the operator
$$
S_{rel}(s) = S_{1}(s) - \psi^{*}S_{2}(s)
$$ 
is trace-class with respect to the Hilbert space completions of
$\FF_{2-s}(\Gamma_1)$ and $\FF_s(\Gamma_1)$, 
then $S_{rel}(s) = 0$ and the diffeomorphism is actually a M\"{o}bius
transformation, i.e. the manifolds $M(\Gamma_{i}) = {\bf H}^{3}/\Gamma_{i}$
are isometric.   

Our results show that the size of the operator $S_{rel}(s)$ 
detects how close to being isometric the quotient manifolds are.  
Recall that the deformation theory of Kleinian groups is
based on the notion of quasi-conformal conjugacy (see section 2).
Our main result is:

\medskip

\noindent {\bf Main Theorem:} {\em
Suppose $\Gamma_{i=1,2}$  are convex co-compact, torsion-free 
Klein\-ian groups so that $M(\Gamma_{i=1,2})$ has infinite hyperbolic
volume.  Let
$$\psi: \Omega(\Gamma_{1}) \rightarrow \Omega(\Gamma_{2})$$  
be an
orientation-preserving $C^{\infty}$-diffeomorphism conjugating 
$\Gamma_{1}$  to $\Gamma_{2}$.  Fix $s \in {\bf C}: \Re(s) = 1, s \neq 1$ 
and let $\epsilon > 0$.  There is 
$K(\epsilon) > 1$ so that $\Vert  S_{rel}(s) \Vert  < \epsilon$
implies that $\Gamma_2$ is a $K(\epsilon)$-quasi-conformal
deformation of $\Gamma_1$, where $K(\epsilon) \rightarrow 1$ as $\epsilon 
\rightarrow 0$. 
}

\medskip

\noindent The norm $\Vert  \cdot \Vert $ in the Theorem is the operator
norm for the $L^{2}$ space of sections. 
Independently,  Douady-Earle \cite{DE86}, Reimann \cite{Rei85} and Thurston
\cite{Thu79} have shown that each $K$-quasi-conformal deformation of a
Kleinian group can be extended to an equivariant $\tilde{K}$-quasi-isometry
of ${\bf H}^{3}$, where $\tilde{K} \rightarrow 
1$ as $K \rightarrow 1$.  Thus the Main Theorem 
says that if $S_{rel}(s)$ is small in the operator norm,
then $M(\Gamma_{1})$ is \lq\lq nearly isometric\rq\rq\ to $M(\Gamma_{2})$.
In particular, $S_{rel}(s) = 0$ implies that $M(\Gamma_1)$ is isometric to
$M(\Gamma_2)$ (a fact which was contained in Perry's result \cite{Per95}).

The plan for this paper is as follows:  Section 2 and  and Section 3
discuss respectively the basics of Kleinian group theory and scattering
theory we will be using. Section 4 contains the proof of the Main Theorem,
as well as various remarks and conjectures.

\smallskip

\noindent{\bf Acknowledgements:}  The authors would like to thank Peter
Perry for explaining the contents of \cite{Per95} to the third author while
Professor Perry was visiting the SUNY-Stony Brook during the Spring of
1993.   We are also indebted to  Richard Canary of the University of
Michigan for helpful and enjoyable conversations on the contents of this
paper.

\section{Kleinian Group Basics}

We will work interchangeably in both the upper half space model and the
ball model of hyperbolic $3$-space; both are denoted by ${\bf H}^{3}$ and  
the distance between $x,y \in {\bf H}^{3}$ is given by $\rho(x,y)$.
Let  $Isom_{+}({\bf H}^{3})$ denote the group of orientation-preserving
isometries of hyperbolic $3$-space ${\bf H}^{3}$ endowed with the
compact-open topology.  A {\em  Kleinian group} $\Gamma$ is a discrete
subgroup of $Isom_{+}({\bf H}^{3})$.   Recall that $Isom_{+}({\bf H}^{3})$
has a natural identification with the   
space of M\"{o}bius transformations $M\ddot{o}b(2)$ 
$$\{ g(z) = \frac{az +b}{cz +d}: \: a,b,c,d \in {\bf C}, ad-bc = 1\}.$$  
Thus $\Gamma$ acts on the
 Riemann sphere $\hat{{\bf C}} = {\bf C} \cup \{\infty\}$ as  
a group of conformal homeomorphisms.  
This action partitions 
$\hat{{\bf C}}$ into two disjoint sets: the limit set and
 the regular set.
Suppose that there is a point $z \in \hat{{\bf C}}$ and a neighborhood
$U$ of $z$, so that
 $\gamma(U) \cap U = \emptyset$ for all $\gamma \in \Gamma - \{id\}$.
 The {\em regular set} of $\Gamma$, denoted by $\Omega(\Gamma)$, is
the (possibly empty) maximal collection of these points.
The {\em limit set} $L_{\Gamma}$ is the complement of the regular set in
$\hat{{\bf C}}$.   The reader is 
referred to \cite{Mas87} for a discussion of the fundamentals in the theory
of Kleinian groups.

Every Kleinian group $\Gamma$ acts discontinuously on ${\bf H}^{3}$, and
there is a natural
geometric model for this action.  Choose a point
$0 \in {\bf H}^{3}$ not fixed by any non-trivial element of $\Gamma$.  The 
{\em Dirichlet polyhedron based at $0$} is the set
$$ P_{0}(\Gamma) = \{x \in {\bf H}^{3}: \: \rho(x,0) \leq \rho(x,\gamma(0))
\; \forall \gamma \in \Gamma\}.$$ 
The intersection of the Euclidean closure of $P_{0}(\Gamma)$ with $\partial
{\bf H}^{3} = \hat{{\bf C}}$ 
is a {\em fundamental domain} ${\cal D}_{\Gamma}$ for the action of
$\Gamma$.  We can form the quotient 
$\Omega(\Gamma)/\Gamma$.   If $\Gamma$ is finitely generated, then 
by the Ahlfors Finiteness Theorem \cite{Ahl64} we can conclude that 
 $\Omega(\Gamma)/\Gamma$  consists of a finite collection of Riemann
surfaces,  each of finite type (i.e. each has finite genus and at most a
finite number of punctures).

A Kleinian group is {\em geometrically finite} if its action on ${\bf
H}^{3}$ admits a finite-sided Dirichlet polyhedron.  We define  $\Gamma$ to
be {\em convex co-compact} if $P_{0}(\Gamma)$ is 
finite-sided, and the Euclidean closure of $P_{0}(\Gamma) \cap \hat{{\bf
C}}$ is bounded away (in the chordal metric) from $L_{\Gamma}$.  
Specifically, in the case that $\Gamma$ is convex co-compact, torsion-free
and $\Omega(\Gamma) \neq \emptyset$,  then $ \overline{M(\Gamma)} = ({\bf
H}^{3} \cup \Omega(\Gamma))/\Gamma$ is a compact $3$-manifold.  The
interior of the $\overline{M(\Gamma)}$ has a complete infinite volume
hyperbolic structure, and the conformal boundary at infinity
$\Omega(\Gamma)/\Gamma$ consists of a finite collection of  compact Riemann
surfaces.  The compactness of $\Omega(\Gamma)/\Gamma$ is a key assumption
in this paper. 
  
Geometrically finite Kleinian groups with non-empty regular sets admit a
space of deformations of at least 1 complex dimension 
via the theory of quasi-conformal mappings (see \cite{Bers74}).   Suppose
$f : \hat{{\bf C}} \rightarrow \hat{{\bf C}}$  
is a quasi-conformal automorphism.   The 
{\em Beltrami coefficient} $\mu$ of $f$ is an element in the unit open ball
in the complex Banach space $L_{\infty}({\bf C})$ of equivalence classes of
bounded measurable functions defined by the condition 
\begin{equation}
\frac{\partial f}{\partial \overline{z}} = 
\mu \frac{\partial f}{\partial z}  
\end{equation}
a.e. in ${\bf C}$ (here $\frac{\partial f}{\partial \overline{z}}$ and
$\frac{\partial f}{\partial z}$ are generalized derivatives).
The {\em dilatation} of $f$ is defined to be 
$K(f) = \frac{1 + ||\mu||_{\infty}}{1 - ||\mu||_{\infty}}$.  Geometrically
a $K$-qc map has the property that it takes infinetesimal circles to
infinitesimal ellipses, such that the ratio of the major to minor axises of
the ellipses is bounded above by $K$.
A fundamental result of Ahlfors-Bers \cite{AB60} says that for any  
$\mu \in L_{\infty}({\bf C}): \: ||\mu||_{\infty} < 1$ there exists a 
unique quasi-conformal automorphism of
the Riemann sphere $\omega^{\mu}$ fixing the points $\{0,1,\infty\}$, and
every quasi-conformal automorphism $f$ with Beltrami coefficient $\mu$ is
of the form $f = \alpha \circ \omega^{\mu}: \: \alpha \in M\ddot{o}b(2).$

Let $\Gamma$ be a finitely generated Kleinian group.
A quasi-conformal automorphism  $f: \hat{{\bf C}} \rightarrow \hat{{\bf
C}}$ that satisfies $f \Gamma f^{-1} \subset M\ddot{o}b(2)$ 
is called a {\em quasi-conformal deformation of $\Gamma$ }. 
Note that if $w^{\mu}$ induces a quasi-conformal deformation of $\Gamma$
and  $\Gamma^{'} = w^{\mu} \circ \Gamma \circ (w^{\mu})^{-1}$, then it is
easy to check that 
$w^{\mu}: \Omega(\Gamma) \rightarrow \Omega(\Gamma^{'})$, 
and hence $w^{\mu}: L_{\Gamma} \rightarrow L_{\Gamma^{'}}$.

Now let $\psi$ be a diffeomorphism $\Omega(\Gamma_1) \to \Omega(\Gamma_2)$
as in Section 1. To extend the mapping $\psi$ to a 
mapping of the whole Riemann sphere we will use 
a foundational result from the theory of Kleinian groups.   
The following is the version of the Marden Isomorphism Theorem \cite{Mar74} we will 
use for this purpose. 
In the statement, ${\bf \overline{H}}^{3}$ denotes the closed $3$-ball.

\begin{theoremSt}
Suppose $\Gamma_{1}$ and $\Gamma_{2}$ are  Kleinian groups such that $\Gamma_{1}$
is convex co-compact and 
there exists an orientation-preserving diffeomorphism
$$\psi: \Omega(\Gamma_{1}) \rightarrow \Omega(\Gamma_{2})$$ 
which induces
an isomorphism $\phi : \Gamma_{1} \rightarrow \Gamma_{2}$.
Then $\psi$ can be extended to a quasiconformal homeomorphism 
$\tilde{\psi}$ of the closed ball ${\bf \overline{H}}^{3}$, also  
inducing $\phi$. 
\end{theoremSt}

\section{Scattering theory}

Scattering theory for the Laplacian can be thought of as a functional
parametrization of the continuous spectrum according to the asymptotic 
behavior of eigenfunctions at infinity.
The interpretation as scattering appears when one translates eigenfunctions
into solutions of either the wave equation or the
Schr\"odinger equation.  Given a choice of incoming solution to one of
these equations (defined asymptotically), the scattering operator provides
the corresponding outgoing solution.  Thus it contains all the information
that can be detected asymptotically on the propagation of waves through the
interior of the manifold.

Scattering theory can be studied on complete Riemannian 
manifolds with certain regular structure at infinity (see \cite{Mel95} 
for details and references).  In our case, the regularity
condition is that in a neighborhood of the conformal boundary at
infinity the metric has the standard form $(dx^2 + dt^2)/t^2$,
$(x,t)\in \R^2\times \R_+$.  This is the reason for the
requirement that  $\Gamma$ be convex co-compact, which disallows cusps.
We hope to extend our results to hyperbolic 3-manifolds with cusps in the 
future.   

The limiting behavior at infinity of eigenfunctions of the Laplacian on
$M(\Gamma) = {\bf H}^3/\Gamma$ is described by sections of complex line
bundles on the conformal boundary.
In particular, different line bundles over $\Omega(\Gamma)/\Gamma$ describe
incoming and outgoing solutions of the wave equation.
We define these line bundles by their spaces of sections, using automorphic
forms.  An {\em automorphic form of weight $s\in\C$} 
on $\Omega(\Gamma)$ is a function $f\in C^\infty(\Omega(\Gamma))$ such that
$$
f(\gamma(x)) |\gamma'(x)|^s = f(x),\qquad\forall \gamma\in\Gamma.
$$
Here $|\gamma'(x)|$ is the {\em conformal dilation} of the
M\"{o}bius transformation $\gamma$ at $x \in \hat{{\bf C}}$.  
Denote by $\FF_s(\Gamma)$ the space of all such forms.
Each $\FF_s(\Gamma)$ corresponds to the space of smooth sections of a
complex line bundle over $\Omega(\Gamma)/\Gamma$.

The scattering operator $S(s)$ 
maps $\FF_{2-s}(\Gamma)$ to $\FF_s(\Gamma)$.  
Its Schwartz kernel can be written as a series,
\begin{equation}
S(s; x,y) = \sum_{\gamma\in\Gamma} {|\gamma'(x)|^s \over
|\gamma(x) - y|^{2s}},
\label{sseries}
\end{equation}
which converges at least for $\Re s>2$.  
Note that $S(s; x,\cdot)$ has weight $s$, so that
\begin{equation}
S(s) f(x) = \int_{D_{\Gamma}} S(s; x,y) f(y) dy,
\end{equation}
is well-defined for $f$ of weight $2-s$, 
where $\DD_{\Gamma}$ is any choice of fundamental domain for the action of
$\Gamma$ on $\Omega(\Gamma)$.
Since $S(s; \cdot,y)$ is of  weight $s$ as well, the scattering operator
maps weight $2-s$ into weight $s$ as claimed.

An automorphic form of weight 2 is
a density, so we have a natural pairing of $\FF_{2-s}(\Gamma)$ 
with $\FF_s(\Gamma)$.  For $\Re s = 1$ complex
conjugatation takes forms of weight $s$ to forms of weight $2-s$, 
so this pairing gives an inner product on $\FF_s(\Gamma)$, $\Re s = 1$:
$$
\pair{f,g} = \int_{\DD_{\Gamma}} \overline{f(x)} g(x) dx.
$$
Accordingly we define the Hilbert space $\HH_\sigma(\Gamma)$ to be the
completion in this inner product of $\FF_{1 + i\sigma}(\Gamma)$.

Of course, the series given in (\ref{sseries}) does not 
necessarily converge for $\Re s = 1$,
the region of interest.  Fortunately, we have the following set of
results.
\begin{theorem}{}
(Mazzeo-Melrose \cite{MM87}, Mandouvalos \cite{Man89}, Perry \cite{Per89}) 

\noindent
Suppose that $\Gamma$ is a convex co-compact co-infinite
volume discrete subgroup of $PSL(2,{\bf C})$.  Let $S(s)$ be the scattering
operator mapping ${\cal F}_{2-s}(\Gamma)$ to ${\cal F}_{s}(\Gamma)$,
defined by (\ref{sseries}) for $\Re s > 2$.
\begin{enumerate}
\item  $S(s)$ has a meromorphic continuation to $s\in\C$, with no poles for
$\Re s = 1$, $s\ne 1$.
\item  $S(s)$ is an elliptic pseudodifferential operator of order
$2-2\Re s$.
\item  In a disk $D\subset \Omega(\Gamma)$, we can write
$$
S(s; x,y) = {1\over |x-y|^{2s}} + \phi(s;x,y), 
\qquad\hbox{\rm for }x,y\in D
$$
where $\phi$ is smooth in $x, y$ and meromorphic in $s$.  For $x$ and $y$
lying in distinct neighborhoods, we have simply that $S(s;x,y)$ is smooth
in $x$ and $y$. 
\end{enumerate}
\end{theorem}

\noindent
{\bf Remarks}
\begin{enumerate}
\item  These results are obtained by studying the scattering kernel
as an asymptotic limit of the resolvent kernel of the Laplacian on
$M(\Gamma)$. After showing in this way that $S(s)$ is an elliptic
pseudodifferential operator for $\Re s \ge 1$, a
parametrix is constructed and used along with a functional relation
$S(2-s) S(s) = I$ to obtain the continuation.
\item  Note that for $\Re s = 1$, $S(s)$ is a zeroth-order operator. 
Since the conformal boundary is compact, $S(1+i\sigma)$ extends to
a bounded operator on $\HH_{-\sigma}(\Gamma)\to
\HH_{\sigma}(\Gamma)$.
\item  The third property, which is crucial for our result, says that
the principal symbol of the scattering operator for $M(\Gamma)$ is the 
same as for ${\bf H}^3$, and in particular independent of $\Gamma$.  
\end{enumerate}

\section{Proof of the Main Theorem}
\newcommand\norm[1]{\Vert#1\Vert}

Recall the setting described in Section 1.
$\Gamma_1$ and $\Gamma_2$ are two convex co-compact, torsion-free
Kleinian groups with non-empty domains of discontinuity
$\Omega(\Gamma_i)$.
For convenience, we'll conjugate $\Gamma_{1}$ so that $0 \in
\Omega(\Gamma_{1})$, and choose a fundamental domain ${\cal
D}_{\Gamma_{1}}$ for the action of $\Gamma_{1}$ on $\Omega(\Gamma_{1})$ so
that $0 \in {\cal D}_{\Gamma_{1}}$.  

Now assume that 
there exists an orientation-preserving diffeomorphism $\psi:
\Omega(\Gamma_1) \to \Omega(\Gamma_2)$ which induces an isomorphism of
$\Gamma_1$ and $\Gamma_2$.   
It is easy to check that under these conditions the diffeomorphism $\psi$
induces a bijection $\psi^*: \FF_s(\Gamma_2) \to \FF_s(\Gamma_1)$, given by
$$
\psi^* f (x) =  (\det D\psi(x))^{s/2}  (f\circ \psi)(x),
$$
where $ \det D(\psi(x))$ is the Jacobian determinant 
of the mapping $\psi$ (see Perry \cite{Per95} section 2).

Let $S_j(s)$, $j=1,2$, denote the scattering operators associated to the
two groups.  We can use $\psi^*$ to define a pullback the scattering
operator $S_2(s)$ to an operator $\psi^* S_2(s): \FF_{2-s}(\Gamma_1) \to
\FF_s(\Gamma_1)$. Thus we can define the {\em relative scattering operator}
$$
S_{rel}(s) = S_1(s) - \psi^*S_2(s),
$$
which we regard by extension as an operator $\HH_{-\sigma}(\Gamma_1) \to
\HH_{\sigma}(\Gamma)$, for $s = 1+i\sigma$.

Our first task is to prove:
\begin{theorem}{}
For $\Re s = 1$, $s\ne 1$, $\epsilon>0$
there exists a number $\delta_s(\epsilon)$ such that
$$
\norm{S_{rel}(s) } < \epsilon \quad\Longrightarrow\quad \psi\hbox{ is }
(1+\delta_s(\epsilon))\hbox{-quasiconformal on } \Omega(\Gamma_{1}),
$$
where $\norm\cdot$ is the operator norm on $\LL(\HH_{-\sigma}(\Gamma_1),
\HH_\sigma(\Gamma_1))$. 
We will also observe that as $\epsilon \downarrow 0$,
$\delta_s(\epsilon) = C\epsilon +  O(\epsilon^2)$, where $C$ is a
constant independent of $s$.
\end{theorem}

\noindent 
{\bf Proof.} 
Fix $s = 1+i\sigma$.
The proof is by contradiction: we assume that $\psi$ is not
$(1+\delta)$-quasiconformal at some point and derive from this a lower
bound (a function of $\delta$) on $\norm{S_{rel}(s)}$.   
We then invert this function to obtain $\delta(\epsilon)$.

Without loss of generality, assume that $\psi$ fails to be
$(1+\delta)$-quasiconformal at $0\in \Omega(\Gamma_1)$.  We can in fact
assume this holds true 
in a neighborhood $|x|<a$ for some $a>0$.  In the course of the proof we
will neglect small $a$ error terms, and the final bound will be shown
independent of $a$.

The principal symbol of $S_1(s)$ is 
$$
a_0(x,\xi) = c_\sigma \;|\xi|^{2i\sigma},
$$
where $c_\sigma$ is a contant with $|c_\sigma|= 1/(2\sigma)$.
The principal symbol of the pullback $\psi^*S_2(s)$ is
$$
a_0'(x,\xi) = c_\sigma \;|A(x)\xi|^{2i\sigma}.
$$
where $A(x)$ is the Jacobian determinant times the inverse of the Jacobian
matrix: 
$$
A(x) = \sqrt{\det D\psi(x)} \cdot (D\psi(x))^{-1}
$$
(note $\det A(x) = 1$).

Noting that $\norm{S_{rel}}^2 = \norm{S_{rel}^*S_{rel}}$, we proceed by
analyzing the oeprator $B = 2\sigma^2 S_{rel}(s)^*S_{rel}(s)$
The principal symbol of $B$ is thus 
$$
b_0(x,\xi) = 1 - \cos \Bigl(\sigma \ln
{|A(x)\xi|^{2}\over |\xi|^{2}}\Bigr).
$$
Let $B_0$ be the pseudodifferential operator with total symbol $b_0(x,
\xi)$, so that $B = B_0 + B_1$ where $B_1$ is of order $-1$.  

Denote by $\lambda(x)$ and $1/\lambda(x)$ be the two eigenvalues of
$A(x)^tA(x)$.  We next quote a standard fact relating the
quasi-conformal factor to the size of $\lambda$ (see \cite{Vai71}, for
example). 

\begin{lemma}{}
With $\lambda(x)$ defined as above, the diffeomorphism
$\psi$ is $K$-quasi\-conformal if and only if $K = \sup\lambda(x)$.
\end{lemma}

As outlined above, we
assume that $\psi$ is not $(1+\delta)$-quasiconformal in a neighborhood
$\{|x|<a\} \subset \Omega(\Gamma_1)$.  We also assume that $\{|x|<a\}
\subset \DD_{\Gamma_1}$.
Our hypothesis for the proof by contradiction thus amounts to 
the assumption that $\lambda(x) > 1+\delta$ for all $|x| < a$.

The lower bound on $\norm{S_{rel}}$ is obtained from this assumption 
by probing with particular $L^2$ states. 
The first serves to localize near $x = 0$, and is defined in ${\cal
D}_{\Gamma_{1}}$ by 
$$
\phi_1(x) = {C_1\over a} \Theta(a-|x|),
$$
where $\Theta$ is the step function and $C_1$ is a constant chosen so that 
$\norm{\phi_1} = 1$. 

For the second, we would like to choose 
\begin{equation}
\phi_2(x) =  {C_2\over a} e^{-|x|^2/2a^2},  \label{phi2}
\end{equation}
where $\xi_0$ is fixed, so that
\begin{equation}
\hat\phi_2(\xi) = C_2'a e^{-a^2|\xi|^2/2}.  \label{hatphi2}
\end{equation}
But as this does not give a smooth section, we must multiply the $\phi_2$
given in (\ref{phi2}) by a smooth cutoff of
$\phi_2$ near the boundary of ${\cal D}_{\Gamma_{1}}$.  This produces an
error term in the $\hat\phi_2$ given in (\ref{hatphi2}), which is uniformly
$O(a^\infty)$ and decreases rapidly in $\xi$.   Thus we may safely neglect
this error by choosing $a$ sufficiently small.   

Observe that $|\pair{S_{rel}(s)\phi_1, S_{rel}(s)\phi_2}| \le
\norm{S_{rel}(s)}^2$.  Recalling that $B_0 + B_1 = 2\sigma^2 S_{rel}(s)^*
S_{rel}(s)$, this means
\begin{equation}
\norm{S_{rel}(s)}^2 \ge {1\over 2\sigma^2}\pair{\phi_1,
(B_0+B_1)\phi_2)} \label{sreltwo}
\end{equation}
Using the fact that $B_1$ is of order $-1$ it is straightforward
to produce a bound 
$$
\norm{B_1\phi_2} \le Ca,
$$
reflecting the fact that $\hat\phi_2$ spreads out as $a$ becomes small.
Since $a$ can be made arbitrarily small, this combines with (\ref{sreltwo})
to give us
$$
\norm{S_{rel}(s)}^2 \ge {1\over 2\sigma^2} |\pair{\phi_1, B_0\phi_2}| 
$$

For ease of exposition in the following integral formulas, let ${\cal D} =
{\cal D}_{\Gamma_{1}}$.  We proceed from the explicit formula
$$
\pair{\phi_1, B\phi_2} = \int\int_{\cal D}  \phi_1(x) b_0(x,\xi)
e^{ix\cdot\xi} \hat\phi_2(\xi)\>dx\>d\xi
$$
Writing $\xi = r\eta$ where $|\eta| = 1$
$$
\pair{\phi_1, B\phi_2} = C\int\int_{\cal D} \Theta(a-|x|) b_0(x,\eta)
re^{ir(x\cdot\eta)} e^{-a^2r^2/2} \>dx\>dr\>d\eta.
$$
Noting that $|\pair{\phi_1, B\phi_2}| \ge |\Re\pair{\phi_1, B\phi_2}|$,
consider
$$
\Re\pair{\phi_1, B\phi_2} = C\int\int_{\cal D} \Theta(a-|x|) b_0(x,\eta) r
\cos(rx\cdot\eta) e^{-a^2r^2/2} \>dx\>dr\>d\eta.
$$
We extract that $r$ integral:
$$
\int_0^\infty r\cos (rx\cdot\eta) e^{-a^2r^2/2} \>dr =
C{1\over a^2}\int_0^\infty r\cos \bigl(r{x\cdot\eta\over a}\bigr)
e^{-r^2/2}
\>dr.
$$
Since $(x\cdot\eta) \le a$, we have
$$
\int_0^\infty r\cos \bigl(r{x\cdot\eta\over a}\bigr)
e^{-r^2/2} \>dr \ge {1\over 4},
$$
so that
$$
|\pair{\phi_1, B\phi_2}| \ge C{1\over a^2} \int \int_{\cal D} \Theta(a-|x|)
b_0(x,\xi) \> dx\> d\eta.
$$

Now consider the $\eta$ integration.  By an appropriate change
of variables we can write $|A(x)\eta|^2 = \lambda(x)\cos^2\theta +
\lambda(x)^{-1} \sin^2\theta$ and $d\eta= d\theta$.
One can check that as a function of $\lambda>1$, the integral
$$
f_\sigma(\lambda) = \int_0^{2\pi} \Bigl[1 - \cos (\sigma \ln
(\lambda \cos^2 \theta + \lambda^{-1}\sin^2\theta)) \Bigr] \;d\theta,
$$
is monotonically increasing.  We thus have
$$
|\pair{\phi_1, B\phi_2}| \ge Cf_\sigma(1+\delta)
$$
Note that this bound is uniform in $a$, justifying our earlier assumptions.

To summarize, we have shown that if $\psi$ is not
$(1+\delta)$-quasiconformal then $\norm{S_{rel}(s)}^2 >
C\sigma^{-2} f_\sigma(1+\delta)$.  It is easy to see that 
for small $\delta$ (with $\sigma$ fixed), $f_\sigma(1+\delta) = C\sigma^2
\delta^2 + O(\delta^3)$.  This completes the proof of Theorem 4.1.
\rightline{\hfill$\square$}
%\finishproof{Theorem 4.1}

\bigskip

\noindent {\bf Proof of the Main Theorem.} 
We have shown that $\psi$ is a 
$K(\epsilon)$-quasi-conformal mapping defined on $\Omega(\Gamma_{1})$,
where $K(\epsilon) = 1+\delta(\epsilon)$.
By assumption $\psi$ induces an isomorphism $\phi:\Gamma_1 \to \Gamma_2$.  
So by the 
Marden Isomorphism Theorem (Theorem 2.1), we can extend
$\psi$ to a  quasi-conformal homeomorphism  
of the closed ball ${\bf \overline{H}}^{3}$ that induces $\phi$.

Recall that the limit set $L_{\Gamma_{1}}$ has measure zero \cite{Ahl64}, 
so we can assume the Beltrami coefficient of
the extended mapping to be zero on $L_{\Gamma_{1}}$.  Thus
the extension of $\psi$ is a 
$K(\epsilon)$-quasi-conformal deformation taking $\Gamma_{1}$ to
$\Gamma_{2}$. 
\rightline{$\square$}

\bigskip
We can use the
Main Theorem to make the following observation concerning the Hausdorff
dimension of the limit sets $L_{\Gamma_{1}}$ and $L_{\Gamma_{2}}$. 
For a set $E \subset \hat{{\bf C}}$ let $D(E)$ be the {\em Hausdorff
dimension} of $E$.  Work of Gehring-V\"{a}is\"{a}l\"{a} \cite{GV73} and
Astala \cite{Ast94} demonstrates that for 
a set $E \subset \hat{{\bf C}}$, and a $K$-quasi-conformal mapping $f$ with
$E$ in its domain, then the set $f(E)$ has Hausdorff dimension bounded
above and below by 
$$ 
\frac{2 D(E)}{2K + (K-1)D(E)} \leq  D(f(E)) \leq \frac{2KD(E)}{2  +
(K-1)D(E)}.
$$
We recall that the quasi-conformal conjugacy 
$\psi: \hat{{\bf C}} \rightarrow \hat{{\bf C}}$ taking $\Gamma_{1}$ to
$\Gamma_{2}$ has the property that $L_{\Gamma_{2}} = \psi(L_{\Gamma_{1}})$.
Thus we have the following corollary:
\begin{corollarySt}
Suppose $\Gamma_{1}$ and $\Gamma_{2}$ are convex co-compact groups
satisfying the conditions given in the Main Theorem.  Then there exists a
$\nu(\epsilon) > 0$ so that 
$$
|D(L_{\Gamma_{1}}) - D(L_{\Gamma_{2}})| < \nu(\epsilon),
$$
where $\nu(\epsilon) \rightarrow 0$ as $\epsilon \rightarrow 0$.
\end{corollarySt}

\clearpage

\noindent{\bf Remarks:}

\begin{enumerate}

\item It is natural to ask whether a small quasi-conformal deformation
implies that the relative scattering operator is small.  In the discussion
below we assume familiarity with the definition of the Poincar\'{e} series
of a Kleinian group $\Gamma$ and with the exponent of convergence
$\delta(\Gamma)$ of this series; we refer the reader to \cite{Nic89} for
discussion of these matters.  

In \cite{BT}, using entirely
two-dimensional techniques, we have been able to show:

\smallskip

\noindent{\bf Theorem:} {\em
Fix $s \in {\bf C}: \: s= 1 + \sigma i, \sigma \neq 0,$ and
 suppose that $\Gamma_{1}$ is a torsion-free non-elementary convex
 co-compact Kleinian group so that $\delta(\Gamma_{1}) < 1$. 
 Then for any $\epsilon > 0$ and for all dilatations $K$
sufficiently close to $1$, each $K$-quasi-conformal deformation  of
 $\Gamma_{1}$ has the property that 
$$ ||S_{rel}(s) || < \epsilon,$$ 
where $|| \cdot ||$ is the operator
norm on ${\cal L}({\cal H}_{-\sigma}(\Gamma_{1}),{\cal
H}_{\sigma}(\Gamma_{1})).$ 
}

\smallskip

The condition $\delta(\Gamma_{1}) <1$ is
topologically restrictive: under the assumptions
given to $\Gamma_{1}$ in the theorem above, results in \cite{CT94} imply
that $M(\Gamma_{1})$  is the interior of a solid handlebody.
However, we believe that a similar result can be shown for two
convex co-compact groups $\Gamma_{i}$ with no assumptions on the exponents
of convergence of the Poincare series.  This would involve a
$3$-dimensional approach using the asymptotic geometry of $M(\Gamma)$ and
the limiting properties of the integral kernel of the resolvent of the
Laplace operator on $M(\Gamma)$ (see \cite{Per89}).

\smallskip

\item Corollary 4.3 shows that for 
convex co-compact Kleinian groups the size of the relative
scattering operator contains information concerning the distortion of the
Hausdorff dimension of the limit sets.  We are thus motivated to ask: {\em
Can the Hausdorff dimension of the limit set be recovered from scattering
data?}

\end{enumerate}

\end{document}